\documentclass[12pt,notitlepage,twoside]{article}
\usepackage{latexsym}
\usepackage{amsmath}
\usepackage{amsfonts}
\usepackage{amssymb}
\renewcommand{\a }{\alpha }
\renewcommand{\b }{\beta }
\renewcommand{\d}{\delta }

\newcommand{\D }{\Delta }

\newcommand{\e }{\varepsilon }
\newcommand{\g }{\gamma}

\newcommand{\G }{\Gamma }
\renewcommand{\l }{\lambda }

\newcommand{\n }{\nabla }

\newcommand{\Sig }{\Sigma}

\renewcommand{\th }{\theta }
\renewcommand{\o }{\omega }
\renewcommand{\O }{\Omega }

\newcommand{\ov}{\overline}

\newcommand{\be}{\begin{equation}}
\newcommand{\ee}{\end{equation}}
\newenvironment{pf}{\noindent{\bf Proof.}\enspace}{%\rule{2mm}{2mm}
\hfill$\Box$\medskip}
\newenvironment{pfn}[1]{\noindent{\bf Proof of {#1}\enspace}}{%\rule{2mm}{2mm}
\hfill$\Box$\medskip}
\newcommand{\R}{\mathbb{R}}
\newcommand{\Q}{\mathbb{Q}}
\newcommand{\Z}{\mathbb{Z}}

\newcommand{\N}{\mathbb{N}}
\newtheorem{thm}{Theorem}[section]
\newtheorem{pro}[thm]{Proposition}
\newtheorem{lem}[thm]{Lemma}
\newtheorem{rem}[thm]{Remark}
\newtheorem{cor}[thm]{Corollary}
\newtheorem{df}[thm]{Definition}
\numberwithin{equation}{section}

\begin{document}

\author{{\sc   Hichem Chtioui \& Mohameden Ould Ahmedou } 
 \footnote{ Corresponding author, { \tt ahmedou@analysis.mathematik.uni-tuebingen.de}}}

\noindent

\title { \Large \textbf{Conformal metrics of prescribed scalar
    curvature on  $4-$manifolds: \\
The degree zero case}
 }

\date{}

\maketitle

{\bf ABSTRACT.-}
{\footnotesize \noindent
\noindent In this paper, we consider the problem of existence and multiplicity
of conformal  metrics  on a riemannian compact $4-$dimensional manifold $(M^4,g_0)$ with   positive scalar curvature. We prove  new exitence criterium which provides existence results for a dense subset of positive functions and generalizes Bahri-Coron and Chang-Gursky-Yang Euler-Hopf type criterium. Our argument gives estimates on the Morse index of the solutions and  has the advantage to extend known existence results. Moreover it  provides  , for generic $K$  {\it Morse Inequalities at Infinity}, which give a
lower bound on the number of metrics with prescribed scalar curvature in
terms of the topological  contribution of its  {\it critical points at
Infinity} to the difference of topology between the level sets of
the associated Euler-Lagrange functional. }

\bigskip


\section{Introduction and  main results}

Let $(M^4,g_0)$ be a compact $4-$dimensional riemannian manifold with positive scalar curvature $R_{g_{0}}$. Given a $C^2$ function $K$ defined on the manifold, the prescribed scalar curvature problem consists of finding a metric $g$, conformally related to $g_0$, such that te scalar curvature of $(M,g)$ is given by te function $K$. Writing $g = u^2 \, g_0$, this amounts to solve the following nonlinear partial differential equation:
$$
{\bf (P_K)} \quad  L_{g_0} \, = \, K \, u^3, \quad u > 0 \mbox{ in } M^4,
$$
where $L_{g_0}$ denotes the conformal Laplacian operator, defined as:
$$
  L_{g_0} := - \D_g u + \, \frac{1}{6} \,u. 
$$

This
problem  has been subject of intensive studies in the last two decades (see  \sloppy \cite{ambr},
\cite{A}, \cite{aujfa}, \cite{BAU}, \cite{AH}, \cite{Blivre},
\cite{B2}, \sloppy \cite{BC1}, \cite{BCCH}, \cite{BCH}, \sloppy
\cite{CGY}, \cite{csL1}, \cite{csL2},
\cite{csL}, \cite{ES}, \sloppy \cite{H}, \cite{yy1}, \cite{L3},
\cite{matthias}, \cite{SZ}, \sloppy \cite{struwe}
and the references therein).\\

Regarding the existence results of the problem $(SC)$, we recall that on $3-$spheres, an Euler-Hopf type criterium for the function $K$ has been obtained by  A. Bahri and J.M. Coron \cite{BC1}, see also Chang-Gursky-Yang \cite{CGY}. Such a criterium  has been generalized for the $4$-spheres by \cite{BCCH} and on higer dimensional spheres only under a closeness to a constant condition \cite{CY} or a flantess condition on the critical points of the function $K \cite{yy1}$.\\
For higher dimensional spheres ($n \geq 7$), A. Bahri \cite{B2}introduced new invariant and discovered new type of existence results. Some of these results have been generalized in  \cite{BCH}.\\
The main difficulty of this problem comes from the presence of the critical Sobolev expoent, which generates blow up and lack of compactness. Indeed the problem enjoys a varitional structure, however the associated Euler Lagrange functional does not satisfy the {\it Palais Smale condition}.  From the variational viewpoint, it is the  occurence of {\it critical points at Infinity}, that are noncompact orbits of the gradient flow, along which the functional remains bounded and its gradient goes to zero, which prevents the use of variational methods.\\
Among approaches developped  to deal with  this problem, we single out is the blow up analysis of some  subcritical approximation combined with the use of the Leray-Schauder  topological degree, approach  developped by R. Schoen \cite{RS1}, Y.Y. Li \cite{yy1}, \cite{L3}, C.S. Lin and C.C. Chen \cite{csL}, \cite{csL1}, \cite{csL2}, among others. The second one is based on  a carefull study of the critical point at Infinity, though a Morse type reduction and the use of their contribution to the topology of the level sets of the associted Euler-Lagrange functional,  has been initiated by A. Bahri and J.M. Coron \cite{BC} and developped through the works of  A. Bahri, \cite{B2} Ben Ayed, Chen, Chtioui, Hammami, see \cite{BCCH}, \cite{BCH}, Ben Ayed, Ould Ahmedou, \cite{BOA}, among others. Other approaches include  perturbations methodes of Chang-Yang \cite{CY} and Ambrosetti \cite{ambr} and the flow approach of M. Struwe \cite{struwe}. \\
In this paper, we revisit this problem to give new existence as well as multiplicity results, extending previous known ones.\\
To state our results we need to introduce some notations and assumptions.\\
We denote by G(a,.) the Green's function of the conformal Laplacian $L_{g_0}$ with pole at $a$ and by $A_a$ the value of its regular part, evaluated at a. \\
Let $0 < K \in C^2(M^4)$ be a  positive function, defined on the manifold $(M^4,g_0$. We say that the function $K$ satisfies the condition ${ \bf (H_0)}$, if $K$ has only nondegenerate critical points and for each critical point $y$, there holds 
$$
\frac{-\D \, K(y)}{3 K(y)} \, - 2 A_y \, \ne 0.
$$
Denoting $\mathcal{K}$ the set of critical point of $K$, we set
$$
\mathcal{K}^+ := \{ y \in \mathcal{K} ; \,  \frac{-\D \, K(y)}{3 K(y)} \, - 2 A_y \, >  \, 0. \}
$$
To each p-tuple $\tau_p := (y_1,\cdots,y_p) \in (\mathcal{K}^+)^p,$ we associate a Matrix $M(\tau_p) = (M_{ij})$ defined by

\begin{eqnarray}\label{f:Matrix}
M_{ii} \, =  \, \frac{-\D \, K(y_i)}{3 K(y_i)^2} \, - 2 \frac{A_{y_i}}{K(y_i)},  \nonumber \\
M_{ij} \,  =  \, \frac{-2 G(y_i,y_j)}{\sqrt{K(y_i)K(y_j)}} \, \mbox{ for } i \ne j.
\end{eqnarray}
We denote by $\rho(\tau_p)$ the least eigenvalue of $M(\tau_p)$ and we say that a function $K$ satisfies the condition ${\bf (H_1)}$ if for every $\tau_p \in (\mathcal{K}^+)^2$, we have that $\rho(\tau_p) \, \ne 0$.\\
We set
\begin{equation}\label{f:finfty}
\mathcal{F_1} := \{ \tau_p = (y_1,\cdots,y_p) \in (\mathcal{K}+)^p \, ; \, \, \rho(\tau_p) \, > 0 \}
\end{equation}
and define an index 
$
\iota: \mathcal{F_1} \to \Z
$
defined by
$$\iota(\tau_p) := p - 1 \, + \, \sum_{i=1}^p (4 - m(K,y_i)),
$$
where $m(K,y_i)$ denotes the Morse index of $K$ at its critical point $y_i$.\\
Now we state our main result.

\begin{thm}\label{t:main}

Let $0 < K \in C^2(M^4)$ be a  positive function satisfying the conditions $(H_0)$ and $(H_1)$.\\
If there exists $k \in \N$ such that
\begin{enumerate}
\item
$$
\sum_{\tau_p \in \mathcal{F_1} ; \iota(\tau_p) \leq k - 1 } (-1)^{\iota(\tau_p)} \, \neq \, 1,
$$
\item
$$
\forall \tau_p \in \mathcal{F_1}, \,  \iota(\tau_p) \neq k    
$$
\end{enumerate}
Then there exists a solution $w$ to   the problem $(P_K)$ such that:
$$
morse(w) \, \leq k,
$$
where $morse(w)$ is the Morse index of $w$, defined as the dimension of the space of negativity of the linerized operator:
$$
\mathcal{L}_w(\varphi) := L_{g_0} (\varphi) \, - 3 w^2 \varphi.
$$
Moreover for generic $K$, it holds
$$
\# \mathcal{N}_k \, \geq | 1 \, - \, |\sum_{\tau_p \in \mathcal{F_1} ; \iota(\tau_p) \leq k - 1 } (-1)^{\iota(\tau_p)} |,
$$
where  $\mathcal{N}_k$ denotes the set of solutions of $(P_K)$ having their Morse indices less or equal $k$.
\end{thm}
Please observe that, taking  in the above $k$ to be $l_{\#} + 1,$ where $l_{\#} $ is the maximal index over all elements of $\mathcal{F_1}$, the second assumption is trivially satisfied. Therefore in this case, we have the following corollay, which recovers previous existence results, see \cite{L3}, \cite{BCCH},\cite{BOA}.  

\begin{cor}\label{c:cor1}
Let $0 < K \in C^2(M^4)$ be a  positive function satisfying the conditions $(H_0)$ and $(H_1)$.\\
If 

$$
\sum_{\tau_p \in \mathcal{F_1} } (-1)^{\iota(\tau_p)} \, \neq \, 1,
$$
Then the problem $(P_K)$ has at least one solution.\\
Moreover for generic $K$, it holds
$$
\# \mathcal{S} \, \geq | 1 \, - \, \sum_{\tau_p \in \mathcal{F_1}  } (-1)^{\iota(\tau_p)} |,
$$
where  $\mathcal{S}$ denotes the set of solutions of $(P_K)$.
\end{cor}

We point out the the main new contribution of Theorem \ref{t:main} is that we  adresse here the case where the total sum in the above corollary  equals 1, but a partial one is not equal 1. The main issue being the possiblity  to use such an information to prove existence of solution to the problem $(P_K)$. To understand the difficulty in adressing such a case, we give, following YY Li \cite{L3}, a new interpretation of the above counting formula in terms of Leray-Schauder degree. Indeed YY Li proved that, under the assumption of corollary \ref{c:cor1}, there exists $R > 0$ such that  the  all solutions of $(P_K)$ remain, for $\a \in (0,1)$ in 
$$
\O_R := \{  u \in C^{2,\a} ; \frac{1}{R} < u < R, ||u||_{C^{2,\a}} < R \}.
$$
It follows that the Leray Schauder degree $deg(v - L^{-1}(K \, v^3)), \O_R, 0)$ is well defined. Moreover it turns out that:
$$
deg(v - L^{-1}(K \, v^3)), \O_R, 0) = 1 \, - \, \sum_{\tau_p \in \mathcal{F_1}  } (-1)^{\iota(\tau_p)}.
$$
Therefore considering the case where the couting formula in corollary \ref{c:cor1}  equals, amounts to considering zero degree case in the above functional analysis approach.\\
Besides the degree interpretation of the couting formula, another interpretation of the fact that the above sum ist different from one, is that the topological contribution of the {critical points at infinity } to te level sets of the associated Euler-Lagrange functional is not trivial. In view of such an interpretation, the above question can be formulated as follows: what happens if the total contribution 
 is trivial, but some  critical points at infinity induce a difference of topology. Can we still use such a topological information to prove existence of solution ? \\
With respect to the above question,   theorem \ref{t:main} gives a sufficient condition to be able to derive from such a local information,  an existence as well as a multiplicity result together with information on the Morse index of the  obtained  solution. At the end of this paper, see, we give a more general condition. Since this  condition involves the critical points at infinity of the variational problem, we have postponed its statement to this end of the paper. \\
As pointed out above, our result does not  only give existence results, butalso,  under generic conditions, gives  a lower bound on the number of solutions of $(P_K)$. Such a result is reminiscent to the celebrated Morse  Theorem, which states that, the number of critical points of a Morse function defined on a compact manifold, is lower bounded in terms of the topology of the underlying manifold. Our resultat can be seen as some sort of {\it Morse Inequality at Infinity}. Indeed it gives a lower bound on the number of metrics with prescribed curvature in terms of the {\it topology at infinity}.\\
The remainder of this paper is organized as follows. In section 2 we set up the variational problem, its critical points at Infinity are characterized in Section 3. Section 4 is devoted to the proof of the main result theorem \ref{t:main} while we give in Section 5 a more general statement than theorem \ref{t:main}.

\bigskip

\bigskip

\begin{center}

{\bf Acknowledgements}

\end{center}

Part of this work has been written when the second author enjoyed the
hospitality of the  {\it Facult\'e des Sciences de Sfax} and { \it Rutgers University, the state University of New Jersey}. He  would like,
in particular  to acknowledge the excellent working conditions in both institutions.

\section{Variational Structure and the lack of compactness}

In this section we recall the functional setting, its 
variational structure  and its main features. Problem $(P_K)$ has a
variational structure. The Euler-Lagrange functional is
 \be \label{00}
J(u)= \frac {\int_M L_{g_0} u \, u}
{\left(\int_{M}K|u|^4\right)^{ 1/2}}
 \ee
defined on $H^1({M},\R)\setminus \{0\}$ equipped with
the norm
$$
||u||^2=\int_{M}
 L_{g_0} u \, u.
$$
We denote by $\Sigma$ the unit sphere of $H^1(M,\R)$
and we set $\Sigma ^+=\{u \in \Sigma  : \, u\geq 0\}$. The
Palais-Smale condition fails to be satisfied for $J$ on $\Sigma
^+$. In order to characterize the sequences failing the
Palais-Smale condition, we need to introduce some notations.\\
Given $a \in M$, we choose a conformal metric 
$$
g_a := u_a^2 \, g
$$
such that $u_a$ depends smoothly on $a$. Let  $x$ be a conformal normal coordinate centered at $a$ and $\varrho > 0$ uniform independent of $a$ such that $x$ is well defined on $B_{2 \varrho}(a)$.\\
We set
$$
\d_{a,\l} \, := \, c_0 \,  \frac{\l}{ 1 \, + \l^2 |x-a|^2}, \, \, x \in B_{\varrho}(a), \, \l > 0,
$$
\noindent  
where $c_0$ ist chosen such that $\d_{a,\l}$ solves the problem
$$
-\D \d_{a,\l} \, = \,  \d_{a,\l}^3 \mbox{ in } \R^4
$$
\noindent
and
\noindent
$$
\hat{\d}_{a,\l}(x) \, := \, u_a(x) \, \o_a(x) \, \d_{a,\l}(x), 
$$
where  $\o_a$ is a cutoff function such that:

$$
 \o_a(x) \, =  \, 1 \, \mbox{ on }  B_{\varrho}(a), \, \,   \o_a(x) = 0 \mbox{ on } M \setminus B_{2 \varrho}(a)
$$
we define $\varphi_{a,\l}$ to be the solution of

$$
L_{g_0}\varphi_{a,\l} \, = \, 8 \hat{\d}_{a,\l}^3. 
$$

\noindent
Setting
$$
H_{a,\l} := \l (\varphi_{a,\l} - \hat{\d}_{a,\l}),
$$
we have that:
\begin{pro}\label{p:h} $\cite{BCCH}$
For $\l$ large, there exists a constant $C =C(\varrho)$ such that:
$$
|H_{a,\l}|_{L^{\infty}} \, \leq C ; \quad  \l |\frac{\partial H_{a,\l}}{\partial \l}|_{L^{\infty}} \leq C ; \quad  \l^{-1} |\frac{\partial H_{a,\l}}{\partial a}|_{L^{\infty}} \leq C .
$$
Moreover for $\varrho$ small and $\l$ large there holds:

\begin{eqnarray}
H_{a,\l}(a) \to A_a \, \, \mbox{ as }  \l \to \infty \\
H_{a,\l}(x) \to G(a,x) \, \, \mbox{ outside } B_{2 \varrho}(a) \, \mbox{ as }  \l \to \infty,
\end{eqnarray}
where $G(a,x)$ is the Green's function of the conformal
subLaplacian $L_{\th}$ and $A_a$ the value of its regular part
evaluated at $a.$
\end{pro}

We define now the set of potential critical points at infinity associated to the functional $J$.\\
For $\e
>0$ and $p\in \N^*$, let us define
\begin{align*}
V(p,\e )=& \{u\in \Sigma /\exists \, a_i \in M^n,  \l _i
>\e^{-1} , \a _i>0 \mbox{ for } i=1,...,p  \, \, \mbox{ s.t. }\\
  &  ||u-\sum_{i=1}^p\a _i \varphi_i||<\e \, , \,  \big|\frac
{\a _i ^{2}K(a_i)}{\a _j ^{2}K(a_j)}-1\big|<\e ,
\mbox{ and } \e _{ij}<\e \}
\end{align*}
where $\varphi_i=\varphi_{(a_i,\l _i)}$ and $\e _{ij}=(\l _i/\l _j +\l _j/\l
_i + \l _i\l _j d(a_i,a_j)^2)^{-1}$.

\noindent For $w$ a  solution of $(P_K)$ we also define
$V(p,\e,w)$ as
\begin{eqnarray}
\{u\in \Sig/\exists \, \a_0>0 \mbox { s. t. } u-\a_0w\in
V(p,\e)\mbox{ and } |\a_0^2
J(u)^2 \, -1|<\e\}.
\end{eqnarray}
The failure of the Palais-Smale condition can be  described
 as follows.

\begin{pro}\label{p:21}\cite{BOA},\cite{BCCH}
 Let $(u_j)\in \Sig^+$ be a sequence such that $\n
J(u_j)$ tends to zero and $J(u_j)$ is bounded. Then, there exist
an integer $p\in \N^*$, a sequence $\e_j>0$, $\e_j$ tends to zero,
and an extracted subsequece of $u_j$'s, again denoted $u_j$, such that
$u_j\in V(p,\e_j,w)$ where $w$ is zero or a solution of $(P_K)$.
\end{pro}
We consider the following minimization problem for $u\in V(p,\e)$
with $\e$ small
\begin{eqnarray}\label{e:51}
\min_{\a _i>0,\, \l _i>0,\,  a_i\in  \mathbb{S}^n}
\bigg\|u-\sum_{i=1}^p\a _i\varphi _{(a_i,\l_i)} \bigg\|_{H^1}.
\end{eqnarray}
We then have the following 
parametrization of the set $V(p,\e )$.
\begin{pro}\label{p:23}\cite{Blivre}, \cite{BC1}, \cite{BCCH} 
For any $p\in \N^*$, there is $\e _p>0$ such that if $\e <\e _p$
and $u\in V(p,\e )$, the minimization problem \eqref{e:51} has a
unique solution (up to permutation). In particular, we can write
$u\in V(p,\e )$ as follows
$$
u=\sum_{i=1}^p\bar{\a }_i\varphi _{(\bar{a}_i,\bar{\l }_i)}+ v,
$$
where $(\bar{\a }_1,...,\bar{\a
}_p,\bar{a}_1,...,\bar{a}_p,\bar{\l }_1,..., \bar{\l }_p)$ is the
solution of \eqref{e:51} and $v\in H^1(\mathbb{S}^n)$ such that
$$(V_0)\quad\qquad\quad ||v|| \leq \e, \quad
(v,\psi)=0 \mbox{ for } \psi\in \bigcup_{i\leq p,\,  j\leq n}
\bigg\{\varphi_i,\frac{\partial\varphi_i}{\partial\l_i},\frac{\partial\varphi_i}{\partial
(a_i)^j} \bigg\},
$$
where $(a_i)^j$ denotes the $j^{th}$ component of $a_i$ and
$(.,.)$ is the inner scalar associated to the norm $\|.\|$.
\end{pro}
In the following we will say that $v\in (V_0)$ if $v$ satisfies
$(V_0)$.
\begin{pro}\label{p:24} \cite{Blivre} \cite{R}
There exists a $C^1$ map which, to each\\ $(\a_1,..., \a_ p,
a_1,..., a_ p, \l_1,..., \l_ p)$ such that $ \sum_{i=1}^p\a _i\varphi
_{(a_i,\l_i)} \in V(p,\e )$ with small $\e $, associates
$\ov{v}=\ov{v}_{(\a_i,a_i,\l_i )}$ satisfying
$$
J\left(\sum_{i=1}^p\a _i\varphi _{(a_i,\l_i)} +\ov{v}\right)= \min_{ v
\in (V_0)} J\left( \sum_{i=1}^p\a _i\varphi _{(a_i,\l_i)}
 +v\right).
$$
Moreover, there exists $c>0 $ such that the following holds
$$
||\ov{v}||\leq c \left(\sum_{i\leq p}(\frac{|\n
K(a_i)|}{\l_i}+\frac{1}{\l_i^2})+\sum_{k\ne r}\e _{kr} (Log (\e
_{kr}^{-1}))^{1/2}\right).
$$
\end{pro}

\noindent Let $w$ be a  solution of $(P_K)$. The
following proposition defines a parameterization of the set
$V(p,\e,w)$.

\begin{pro}\label{p:w-parametrisierung} \cite{B2}
There is $\e_0>0$ such that if $\e\leq \e_0$ and $u\in V(p,\e,w)$,
then the problem
 $$
 \min_{\a_i>0,\,,\,  \l_i>0,\, \,  a_i\in  {M} , \, \, h\in
T_w(W_u(w))}  \big|\big|u-\sum_{i=1}^p\a_i
\varphi_{(a_i,\l_i)}-\a_0(w+h) \big|\big|
 $$
has a unique solution $(\ov{\a},\ov{\l},\ov{a}, \ov{h} )$. Thus,
we write $u$ as follows:
$$u=\sum_{i=1}^p\ov{\a}_i\varphi_{(\ov{a}_i,\ov{\l}_i)}
+\ov{\a}_0(w+\ov{h})+v,$$ where $v$ belongs to
$H^1(M)\cap T_w(W_s(w))$ and it satisfies $(V_0)$,
$T_w(W_u(w))$ and $T_w(W_s(w))$ are the tangent spaces at $w$ to
the unstable and stable manifolds of $w$.
\end{pro}

\section{Critical points at Infinity of the variational problem}
Following A. Bahri we set the following definitions and
notations
\begin{df}
{\it A critical point at infinity} of $J$ on $\Sig^+$ is a limit
of a flow line $u(s)$ of the equation:
$$
\begin{cases}
\frac{\partial u}{\partial s} = - \n J(u) \\
u(0) = u_0
\end{cases}
$$
such that $u(s)$ remains in $V(p,\e(s),w)$ for $s \geq s_0$.\\
Here $w$ is either zero or a solution of $(P_K)$ and $\e(s)$ is
some function tending to zero when $s \to \infty$. Using
Proposition \ref{p:w-parametrisierung}, $u(s)$ can be written as:
$$
u(s) \, = \,  \sum_{i=1}^p\a_i(s) \, \varphi_{(a_i(s),\l_i(s))} +
\a_0(s)(w + h(s)) \, + v(s).
$$
Denoting $a_i := \lim_{s \to \infty} a_i(s)$ and $\a_i=\lim_{s\to
\infty}\a_i(s)$, we denote by
$$(a_1,\cdots,a_p,w)_\infty \mbox{ or } \sum_{i=1}^p\a_i \,
\varphi_{(a_i,\infty)} + \a_0w$$ such a critical point at infinity. If
$w \ne 0$ it is called of {\it $w$-type}.
\end{df}

\subsection{Ruling out the existence of critical point at Infinity
in $V(p,\e,w)$ for $w \neq 0$ }

The aim of this section is to prove  that, given a funktion $K$ a $C^2$ positive funktion satisfying the condition of theorem \ref{t:main} and $w$  a solution of $(P_K)$. Then for each $p \in \N$, there are no critical point or critical point at infinity of $J$ in the set $V(p,\e,w)$. The reason is that there exists a pseudogradient of $J$ such that the Palais Smale condition is satisfied along the decreassing flow lines.\\
In this section, for $u\in V(p,\e,w)$, using Proposition
\ref{p:w-parametrisierung}, we will write $u=\sum_{i=1}^p\a _i\, \varphi
_{(a_i,\l _i)}+\a_0(w+h)+v$.
\begin{pro}\label{p:32}
For $\e >0$ small enough and $u=\sum_{i=1}^p\a _i\varphi _{(a_i,\l
_i)}+\a_0(w+h)+v\in V(p,\e,w )$, we have the following expansion
\begin{align*}
J(u) = & \frac{S_4\sum_{i=1}^p \a_i^2+\a_0^2||w||^2}{
(S_4\sum_{i=1}^p  \a _i
^4K(a_i)+\a_0^4||w||^2 )^\frac{1}{2}}
\left[1-c_2\a_0\sum_{i=1}^p
\a_i\frac{w(a_i)}{\l_i}\right.\\
& -c_2\sum_{i\ne j} \a _i\a _j\e _{ij}
 +f_1(v)+Q_1(v,v) +f_2(h) +\a_0^2Q_2(h,h)\\
& \left.+o\biggl(\sum_{i\ne j}\e _{ij} +\sum_{i=1}^p\frac{1} {\l
_i }+||v||^2+ ||h||^2\biggr)\right]
\end{align*}
where
\begin{align*}
Q_1(v,v) = & \frac{1}{\gamma_1}||v||^2-\frac{3}{\b_1}
\int_{M^4}K\left(\sum_{i=1}^p(\a _i\varphi
_i)^2+(\a_0w)^2\right)
v^2, \\
Q_2(h,h)= & \frac{1}{\gamma_1}||h||^2-\frac{3}{\b_1}
\int_{M^4}K(\a_0w)^2 h^2,\\
f_1(v) = & -\frac{1}{\b_1 }\int_{M^4}K(\sum_{i=1}^p\a
_i\varphi _i)^3 \, v,  \\
f_2(h)= &
\frac{\a_0}{\gamma_1}\sum_i\a_i(\varphi_i,h)-\frac{\a_0}{\b_1}
\int_{M^n}K(\sum_i\a_i\varphi_i+\a_0 w)^3h,\\
c_2= & c_0^4\int_{\R^4}\frac{dx}{(1+
|x|^2)^{3}} , \qquad  S_4 = c_0^4
\int_{\R^4}\frac{dx}{(1+|x|^2)^4}\\
 \b_1  = & {S_4}(\sum_{i=1}^p
\a _i^4 K(a_i))+\a_0^4||w||^2, \qquad
\gamma_1 = {S_4}(\sum_{i=1}^p \a _i^2)+\a_0^2||w||^2.
\end{align*}
\end{pro}
\begin{pf}
To prove the proposition, we 
need to estimate
$$
N(u)= ||u||^2 \quad \mbox{ and }\quad
D^2=\int_{M^n} K(x) u^4.
$$
We observe first that, expanding $N(u)$, we have that
$$
\sum_{i=1}^p\a_i^2||\varphi_i||^2+2\a_i\a_0<\varphi _i, w+h> +\a_0^2
(||h||^2+||w||^2) +||v||^2 +\sum_{i\ne j}\a_i\a_j<\varphi _i,\varphi_j>.
$$
Now it follows from \cite{BCCH} and elementary computations that
\begin{align}
& ||\varphi_i ||^2= S_4 \, + 2 \o_3 \, \frac{H_{a_i,\l_i}(a_i)}{\l_i^2}  \\ & <\varphi_i,\varphi_j> \, = \,  \frac{2\o_3 \, H_{a_j,\l_j}(a_i)}{\l_i \, \l_j} \, +  \, 
 c_2 \e _{ij}(1 + o(1)) , \quad \mbox{ for } i\ne j, \label{e:o1}\\
 & (\varphi_i,w)=\int_{M^4}w\varphi_i ^{3}=
c_2\frac{w(a_i)}{\l_i}+ o(\frac{1}{\l_i}).
\end{align}
 Therefore
\begin{align}\label{N}
 N=  & \gamma_1+2\a_0\sum_{i=1}^pc_2
\a_i\frac{w(a_i)}{\l_i}+\a_i(\varphi_i,h) +c_2\sum_{i\ne j}\a
_i \a _j\e _{ij}\\
 &  +\a_0^2||h||^2+ ||v||^2+
o\biggl(\sum_{i=1}^p\frac{1}{\l_i}+ \sum_{i\ne j}\e
_{ij}\biggr).\notag
\end{align}

Now concerning the denominator, we compute it as follows
\begin{align}
D^2=&\int K (\sum_{i=1}^p\a _i\varphi _i
)^4+\int K
(\a_0w)^4\\
& + 4 \a_0 \, \int K (\sum_{i=1}^p\a _i\varphi
_i)^4 w + 4 \a_0^3\int K
(\sum_{i=1}^p\a _i\varphi _i)w^3\notag\\
&  + 4 \int K
(\sum_{i=1}^p\a _i\varphi _i+\a_0w)^4(\a_0h+v)  \notag\\
 & + 6\int K(\sum_{i=1}^p\a_i \varphi _i+\a_0w)^2
(\a_0^2h^2+v^2+2\a_0hv)\notag\\
 & +O(\sum \int
 w^2\varphi_i^2+w^2\varphi_i^2)+O(||v||^{3}+\|h\|^3).\notag
\end{align}
 Observe that
\begin{align}
\int_{M^4}K (\sum_{i=1}^p\a _i\varphi _i)^4 = &
\sum_{i=1}^p\a_i ^4K(a_i) S_4\\
& + 4   c_2\sum_{i\ne j} \a_i^4\a_j
K(a_i)\e_{ij} +O( \frac{1}{\l _i^2})+o(\e _{ij}), \notag
\end{align}

\begin{align}
& \int_{M^4}K w^4=||w||^2; \quad
\int_{M^4} Kw^3 \, \d_i =
c_2\frac{w(a_i)}{\l_i }+  ,\\
& \int_{M^4} K(\sum \a_i\varphi_i) ^3 \, w=c_2\sum
\a_i^3
K(a_i)\frac{w(a_i)}{\l_i } \, + \, o(\frac{1}{\l_i}),\\
& \int_{M^4}\varphi_i^2 w^2+\varphi_i^2 w^2
=o(\frac{1}{\l_i }),
\end{align}

\begin{align}
\int_{M^4} K(\sum \a_i \varphi_i +\a_0 w)^2 \,  v \, h & =
O\left(\int
\left(\sum \varphi_i^2 \,  +w^{-1}\sum  \varphi_i\right)|v||h|\right)\notag\\
& =O\left(\|v\|^3+\|h\|^3+1/\l_i^{3}\right),
\end{align}
where we have used that $v\in T_w(W_s(w))$ and $h$ belongs to
$T_w(W_u(w))$ which is a finite dimensional space. Hence it
implies that $\|h\|_\infty \leq c\|h\|$.\\
Concerning the linear form in  $v$, since $v\in T_w(W_s(w))$, it
can be written as
\begin{align}
\int_{M^4}K (\sum_{i=1}^p & \a _i \varphi _i +\a_0
w)^{3}v\notag\\ & =\int K (\sum_{i=1}^p\a _i\varphi
_i)^{3}v
+O\biggl(\sum_{i=1}^p \int(\varphi_i^2 \, w +\d_i w ^2 )|v|\biggr)\notag\\
 & =f_1(v)+O\left(\frac{||v||}{\l_i}\right).
\end{align}
Finally, we have
\begin{align}
 \int K(\sum_{i=1}^p\a _i\varphi _i+\a_0w)^2 \, h^2= &
\a_0^2 \int
Kw^2 \, h^2 + o(||h||^2)\\
 \int K(\sum_{i=1}^p\a _i\varphi_i+\a_0w)^2 \, v^2= &
\sum_{i=1}^p\int K(\a _i\varphi_i)^2 \, v^2+\a_0^2
\int Kw^2 v^2\notag\\
 & +o(||v||^2).\label{v2}
\end{align}
Combining \eqref{N},...,\eqref{v2}, the result follows.
\end{pf}

Now, we state  the following lemma which is proved for the
dimensions $n\geq 7$ in \cite{B2} in the case of the spheres but the proof works virtually in our case.
\begin{lem} We have \\
(a) $Q_1(v,v)$ is a quadratic form positive definite in\\
\centerline{ $E_v=\{v\in H^1(M^4)/v\in T_w(W_s(w)) \mbox{
and } v \mbox{ satisfies } (V_0)\}$.}
 (b) $Q_2(h,h)$ is a quadratic form
negative definite in $T_w(W_u(w))$.
\end{lem}

\begin{cor}\label{c:} \cite{B2} Let $u=\sum_{i=1}^p\a _i\d _{(a_i,\l
_i)}+\a_0(w+h)+v\in V(p,\e,w )$. There is an optimal
$(\ov{v},\ov{h})$ and a change of variables $v-\ov{v}\to V$ and
$h-\ov{h}\to H$ such that %$J$ reads as
$$J(u) = J\left(\sum_{i=1}^p\a _i\d _{(a_i,\l
_i)}+\a_0 w+ \ov{h}+\ov{v}\right)+||V||^2-||H||^2.
$$
Furthermore we have the following estimates
$$
||\ov{h}|| \leq \sum_i\frac{c}{\l_i }\, \, \,
\mbox{ and }\, \, \,  ||\ov{v}||\leq c\sum_i\frac{|\n
K(a_i)|}{\l_i}+\frac{c}{\l_i^2}+c
\sum\e_{kr}(Log\e_{kr}^{-1})^{\frac{1}{2}},
$$
\begin{align*}
J(u)= & \frac{S_n\sum_{i=1}^p\a_i^2+\a_0^2||w||^2}{
(S_n\sum_{i=1}^p\a _i ^4 \, K(a_i)+\a_0^4
||w||^2 )^{\frac{1}{2}}} \left[1-c_2\a_0\sum_{i=1}^p\a_i
\frac{w(a_i)}{\l_i }\right.\\
& -c_2\sum_{i\ne j}\a _i\a _j\e _{ij} \left.+o\left(\sum_{i\ne
j}\e _{i j}+  \sum_{i=1}^p\frac{1} {\l _i
}\right)\right]  +||V||^2-||H||^2.
\end{align*}
\end{cor}
\begin{pf}
The expansion of $J$ with respect to $h$ (respectively to $v$) is
very close, up to a multiplicative constant, to $Q_2(h,h)+f_2(h)$
(respectively $Q_1(v,v)+f_1(v)$). Since $Q_2$ is negative definite
(respectively $Q_1$ is positive definite), there is a unique
maximum $\ov{h}$ in the space of $h$'s (respectively a unique
minimum $\ov{v}$ in the space of $v$). Furthermore, it is easy to
derive $||\ov{h}||\leq c||f_2||$ and $\|\ov{v}\|\leq c\|f_1\|$.
The estimate of $\ov{v}$ follows from Proposition \ref{p:24}. For
the estimate of $\ov{h}$, we use the fact that for each $h\in
T_w(W_u(w))$ which is a finite dimensional space, we have
$\|h\|_\infty \leq c \|h\|$. Therefore, we derive that
$\|f_2\|=O(\sum \l_i^{-1})$. Then our result follows.
\end{pf}

\noindent
Now we state the following corollary, which follows   immediately from the above corollary and
the fact that $w > 0$ in ${M^4}$.

\begin{cor}\label{c::}
Let $K$ be a $C^2$ positive function and let $w$ be a
nondegenerate critical point of $J$ in $\Sig^+$. Then, for
each $p\in \N^*$, there is no critical
points or critical points at infinity in the set $V(p,\e,w)$, that
means we can construct a pseudogradient of $J$ so that the
Palais-Smale condition is satisfied along the decreasing flow
lines.
\end{cor}

\noindent
Now once mixed critical points at Infinity is ruled out,  it follows from \cite{BCCH} and \cite{BOA}, that the critical points at infiity are in one to one correspondence with the elements of the set $\mathcal{F_1}$  defined in \eqref{f:finfty}. that is  a critical point at infinity corresponds to $\tau_p:= (y_1,\cdots,y_p) \in (\mathcal{K}^+)^p$ such that the realted Matrix $M(\tau_p)$ defined in \eqref{f:Matrix} is positive definite. Such a critical point at infinity will be denoted by $\tau_p^{\infty} := (y_1,\cdots,y_p)_{\infty}.$ \\
Like a usual critical point, it is associated to  
 a {\it critical point at infinity} $x_{\infty}$ of the problem $(P_K)$, which are combination of classical critical points with a $1-$dimensional assymptote,   stable and unstable
manifolds, $W_s^{\infty}(x_{\infty})$ and $W_u^{\infty}(x_{\infty})$.   These manifolds can be easily described once a Morse type
reduction is performed, see \cite{B2}, \cite{BCCH}. The stable amnifold is, as usual, defined to be the set of points attracted by the asymptote. The unstable one is a shadow object, which is the limit of  $W_u(x_{\l})$, $x_{\l}$ being the critical point of the reduced problem and $W_u(x_{\l})$ its associated unstable manifolds. Indeed the flow in this case splits the variable $\l$ from the other variables near $x_{\infty}$.\\
In the following defiition, we extend the notation of domination of critical points to critical points at Infinity.
\begin{df}
$z_\infty$ is said to be dominated by another critical point at
infinity $z'_\infty$ if
$$W_u(z'_\infty)\cap W_s(z_\infty)\ne \emptyset.$$
\end{df}
If we assume that the intersection is transverse, then we obtain
$$index(z'_\infty)\geq index(z_\infty)+1.$$

\section{Proof of the main result}

This section is devoted to the proof of the main result of this paper, theorem \ref{t:main}.
\begin{pfn}{\bf Theorem \ref{t:main} }

Setting 
$$
l_{\#} := \sup \{ \iota(\tau_p); \, \tau_p \in \mathcal{F_1}  \}
$$
For $l \in \{0,\cdots,l_{\#} \}$ we define the following sets:
\begin{equation}\label{f:xl}
X^{\infty}_l := \cup_{\tau_p \in \mathcal{F_1} ; \, \iota(\tau_p) \leq l} \ov{W^{\infty}_s(\tau_p^{\infty})},
\end{equation}
where $W_u^{\infty}(\tau_p^{\infty})$ is the unstable manifold associated to the critical point at infinity $\tau_p^{\infty}$.
and
\begin{equation}\label{f:cl}
C(X^{\infty}_l) := \{ t \, u \, + \, (1 -t) \, (y_0)_{\infty},\,  t \in [0,1],\,  u \in    X^{\infty}_l \},
\end{equation}
where $y_0$ is  a global maximun of $K$ on the manifold $M^4$.\\
By a theorem of Bahri-Rabonowitz \cite{BR}, it follows that:

$$
\ov{W^{\infty}_s(\tau_p^{\infty})} \, = \, W^{\infty}_s(\tau_p^{\infty}) \,  \cup \,  \cup_{x_{\infty} < \tau_p^{\infty}} W^{\infty}_s(x_{\infty}) \, \cup \,   \cup_{w < \tau_p^{\infty}} W_u(w),
$$
where $x_{\infty}$ is a critical point at infinity dominated by $\tau_p^{\infty}$ and $w$ is a solution of $P_K$ dominated by $\tau_p^{\infty}.$ By tranversality arguments for we assume that the index of $x_{\infty}$ and the Morse index of $w$ are no bigger than $l$. Hence
$$
X^{\infty}_l \, = \, \cup_{\iota(\tau_p) \leq l} W^{\infty}_s(\tau_p^{\infty}) \,\cup_{w < \tau_p^{\infty}} W_u(w).
$$
It follows  that $X^{\infty}_k$ is a stratified set of top dimension $ \leq l$. Without loss of generality we may assume it  equals to $l$ therefore $C(X^{\infty}_k)$ is  also a stratified set of top  dimension $l + 1$.\\
Now we use the gradient flow of $- \n J$ to deform $C(X^{\infty}_k)$. By tranversality arguments we can assume that the deformation avoids all critical as well as critical points at Infinity having their Morse indices greater than $l + 2$. It follows then by a a Theorem of Bahri and Rabinowitz \cite{BR}, that  $C(X^{\infty}_k)$ retracts by deformation on the  set 
\begin{equation}\label{f:fo}
U :=  X^{\infty}_l \,  \cup \, \cup_{\iota(x{_{\infty}}) = l + 1} W^{\infty}_u(x_{\infty})  \, \cup \, \cup_{w < \tau_p^{\infty}} W_u(w). 
\end{equation}

\noindent
Now taking $l = k - 1$ and using that by assumption of theorem \ref{t:main}, there are no critical pointa t infinity with index $k$, we derive that
$C(X^{\infty}_k)$ retracts by deformation onto
\begin{equation}\label{f:rpdef}
Z_k^{\infty} := X^{\infty}_k  \, \cup  \cup_{w ; \n J(w) = 0 ; w \mbox{ dominated by } C(X^{\infty}_k)  } W_u(y).
\end{equation}
Now observe that, it follows from the above deformation retract that the problem $(P_K)$ has necessary a solution $w$ with $m(w) \leq k$. Otherwise  it follows from \eqref{f:rpdef}that 

$$
1 \, = \chi(Z_k^{\infty}) \, = \,  \sum_{\tau_p \in \mathcal{F_1} \, ; \iota(\tau_p) \leq  k - 1} (-1)^{ \iota(\tau_p)} \, ,
$$ 
where  $\chi$ denotes the Euler Characteristic. Such an equality contradicts the assumption 2 of the theorem. \\
Now for generic $K$, it follows from the Sard-Smale Theorem that all solutions of $(P_K)$ are nondegenerate solutions, in the sens that their associated linearized operator does not admit zero as an eingenvalue. See \cite{SZ}.\\
We derive now from \eqref{f:rpdef}, taking the Euler Characteristic of both sides that:
$$
1 \, = \chi(Z_k^{\infty}) \, = \,  \sum_{\tau_p \in \mathcal{F_1} \, ; \iota(\tau_p) \leq  k - 1} (-1)^{ \iota(\tau_p)} \,  + \, \sum_{w < X_k^{\infty} ; \n J(w) = 0} (-1)^{m(w)}.
$$ 
It follows then that
$$
| 1 -\sum_{\tau_p \in \mathcal{F_1} \, ; \iota(\tau_p) \leq  k - 1} (-1)^{ \iota(\tau_p)} | \leq  \sum_{w  ; \n J(w) = 0, m(w) \leq k} (-1)^{m(w)} \leq \mathcal{N}_k,
$$
where $\mathcal{N}_k $ denotes the set of solutions of $(P_K)$ having their morse indices $\leq k$. 
\end{pfn}

\section{A general existence result}

In this last section of this paper, we give a generalization of theorem \ref{t:main}. Namely instead of assuming that there are no   critical point at infinity of index $k$, we assume that the interesction number modulo 2, between the suspension of the complex at infinity of order $k$, $C(X_k^{\infty})$ and the stable manifold of all critical points at infinity of index $k + 1$ is equal zero. More precisely, for $\tau_p \in \mathcal{F_1}$ such that $\iota(\tau_p) \, = k $, we define the following intersection number:
$$
\mu_k(\tau_p) := C(X_{k - 1}^{\infty}) \, . \, W_s^{\infty}(\tau_p^{\infty}) \, (\mbox{mod} 2).
$$
Observe that this intersection number is well defined since we may assume by transversality that:
$$
\partial C(X_k^{\infty}) \, \cap \, W_s^{\infty}(\tau_p^{\infty}) \, = \, \emptyset. 
$$
indeed
$$
dim (\partial C(X_k^{\infty})) \, = \,  k - 1, \mbox{ while } dim (W_s^{\infty}(\tau_p^{\infty})) = 4 - k.
$$ 
We are now ready to state the following existence result:

\begin{thm}\label{t:a}

Let $0 < K \in C^2(M^4)$ be a  positive function satisfying the conditions $(H_0)$ and $(H_1)$.\\
If there exists $k \in \N$ such that
\begin{enumerate}
\item
$$
\sum_{\tau_p \in \mathcal{F_1} ; \iota(\tau_p) \leq k - 1 } (-1)^{\iota(\tau_p)} \, \neq \, 1,
$$
\item
$$
\forall \tau_p \in \mathcal{F_1}, \mbox{ such that } \iota(\tau_p) = k, \mbox{ there holds } \mu_k(\tau_p) \, = \, 0.    
$$
\end{enumerate}
Then there exists a solution $w$ of  the problem $(P_K)$ such that:
$$
morse(w) \, \leq k,
$$
where $morse(w)$ is the Morse index of $w$.\\
Moreover for generic $K$, it holds
$$
\# \mathcal{N}_k \, \geq | 1 \, - \, |\sum_{\tau_p \in \mathcal{F_1} ; \iota(\tau_p) \leq k - 1 } (-1)^{\iota(\tau_p)} |,
$$
where  $\mathcal{N}_k$ denotes the set of solution of $(P_K)$ having their Morse indices less or equal $k$.
\end{thm}

\begin{pf}
THe proof goes  along with the proof of theorem \ref{t:main},therefore  we will only sketch the differences. Keeping the notation of the proof of theorem \ref{t:main}, we observe that, since 

$$
\forall \tau_p \in \mathcal{F_1}, \mbox{ such that } \iota(\tau_p) = k, \mbox{ there holds } \mu_k(\tau_p) \, = \, 0,
$$
we may assume that the deformation of $C_k^{\infty}$ along any pseudogradient flow of $-J$, avoids all critical points at infinity having their Morse indices equal to $k$. It follows then from \eqref{f:fo}
that 
$C(X^{\infty}_k)$ retracts by deformation onto
\begin{equation}\label{f:rpdefa}
Z_k^{\infty} := X^{\infty}_k  \, \cup  \cup_{w ; \n J(w) = 0 ; w \mbox{ dominated by } C(X^{\infty}_k)  } W_u(y).
\end{equation}
Now the remainder of the proof is identical to the proof of theorem \ref{t:main}.
\end{pf}

{\small

\bigskip
\bigskip
\newpage

\centerline{\sc Hichem Chtioui}
\centerline {D\'epartement de Math\'ematiques}
\centerline{ Facult{\'e} des Sciences de Sfax}
\centerline{ Route Soukra, Sfax, Tunisia.}
\centerline{\tt Email: hichemchtioui2003@yahoo.fr}
\bigskip
\centerline{\sc and }
\bigskip

\centerline{ \sc Mohameden  Ould Ahmedou}
 \centerline{Mathematisches institut}
\centerline{Universit\"at T\"ubingen}
\centerline{ Auf der Morgenstelle 10}
\centerline{ D-72076 Tubingen, Germany.}
\centerline{\tt Email:ahmedou@analysis.mathematik.uni-tuebingen.de }
}

\end{document}